\documentclass[12pt]{amsart}

\usepackage{color}
\usepackage[centertags]{amsmath}
\usepackage[margin=1in]{geometry}
\usepackage{amsfonts}
\usepackage{amsthm}
\usepackage{newlfont}
\usepackage{verbatim}
\usepackage{amscd}
\usepackage{amsgen}
\usepackage{youngtab}
\usepackage{amssymb}
\usepackage{stmaryrd}
\usepackage{amssymb,amsmath}
\usepackage{amscd}
\usepackage{enumerate}
\usepackage{color}
\usepackage{xypic}
\usepackage{graphicx}

\newlength{\defbaselineskip} 
\theoremstyle{plain}
\newtheorem{thm}{Theorem}[section]
\newtheorem{cor}[thm]{Corollary}
\newtheorem{con}[thm]{Conjecture}
\newtheorem{df}[thm]{Definition}
\newtheorem{lema}[thm]{Lemma}
\newtheorem{obs}[thm]{Proposition}
\newtheorem{exm}[thm]{Example}

\newtheorem*{tthm}{Main Theorem}

\newtheorem{rem}[thm]{Remark}
\newtheorem{pr}{Algorithm}
\newtheorem{fact}[thm]{Fact}
\theoremstyle{definition} 
\theoremstyle{definition} \newtheorem{ex}{Example}[section] %

 \numberwithin{equation}{section}
\def\p{\mathbb{P}}
\def\P{\mathbb{P}}

\def\n{\mathbb{N}}
\def\Z{\mathbb{Z}}
\def\N{\mathbb{N}}
\def\c{\mathbb{C}}
\def\C{\mathbb{C}}

 \DeclareMathOperator{\h}{ht}

\def\p{\mathbb{P}}
\def\a{\mathbb{A}}
\def\ob{\begin{obs}}
\def\kob{\end{obs}}
\def\dow{\begin{proof}}
\def\kdow{\end{proof}}
\def\kwadrat{\hfill$\square$}
\def\tw{\begin{thm}}
\def\ktw{\end{thm}}
\def\hip{\begin{con}}
\def\khip{\end{con}}
\def\lem{\begin{lema}}
\def\klem{\end{lema}}
\def\ex{\begin{exm}}
\def\prog{\begin{pr}}
\def\kprog{\end{pr}}
\def\wn{\begin{cor}}
\def\kwn{\end{cor}}
\def\uwa{\begin{rem}}
\def\kuwa{\end{rem}}
\def\kex{\end{exm}}
\def\dfi{\begin{df}}
\def\kdfi{\end{df}}
\def\g{\mathfrak{g}}\def\h{\mathfrak{h}}
\def\n{\mathfrak{n}}
\def\pp{\mathfrak{p}}
\def\b{\mathfrak{b}}
\definecolor{zielony}{rgb}{0.5, 0.9, 0.1}
\definecolor{czerwony}{rgb}{0.9, 0.2, 0.1}
\definecolor{niebieski}{rgb}{0.3, 0.1, 0.9}

\xyoption{line}
\def\fa{\begin{fact}}
\def\kfa{\end{fact}}
\DeclareMathOperator{\sgn}{sgn}
\author{Laurent Manivel}
\address{Institut Fourier, Universit\'e Joseph Fourier, 100 rue des Math\'ematiques,
F-38400 Saint Martin d'H\`eres, France}
\email{Laurent.Manivel@ujf-grenoble.fr}
\author{Mateusz Micha\l ek}
\address{
Max Planck Institute for Mathematics,
Vivatsgasse 7,
53111 Bonn,
Germany
\linebreak
Mathematical Institute of the Polish Academy of Sciences, \'{S}niadeckich 8, 00-956 Warszawa, Poland}
\email{wajcha2@poczta.onet.pl}
\thanks{The first author is supported by the Narodowe Centrum Nauki grant UMO-2011/01/N/ST1/05424 "Representation theory and secants of homogeneous varieties".}

\title{Secants of minuscule and cominuscule minimal orbits}
\begin{document}
\maketitle
\begin{abstract}
We study the geometry of the secant and tangential variety of a cominuscule and minuscule variety, e.g. a Grassmannian or a spinor variety.
Using methods inspired by statistics we provide an explicit local isomorphism with a product of an affine space with a variety which is the Zariski
closure of the image of a map defined by generalized determinants. In particular, equations of the secant or tangential
variety correspond to relations among generalized determinants. We also provide a representation theoretic decomposition of cubics in the ideal of
the secant variety of any Grassmannian.
\end{abstract}

\section{Introduction}
The aim of the article is to investigate the properties of the secant variety of the minimal orbit in a minuscule and cominuscule representation
of a semi-simple complex Lie group.
The prototypical examples of such varieties are the Grassmannians. The Grassmannian of $k$ dimensional subspaces of an $n$ dimensional vector
space $V$ is the image of the map
 $$  \{\text{nondegenerate }k\times n\text{ matrices}\} \rightarrow \P(\bigwedge^kV) $$
$$  M \mapsto  [ \text{all maximal minors of }M ] .$$

 \medskip
Moreover, we can parameterize an affine open chart of the Grassmannian by
$$\{ k\times (n-k)\text{ matrices}\}\rightarrow \a^{{n\choose k}-1}\subset\P(\bigwedge^kV)$$
$$M\mapsto ( \text{all minors of }M ).$$

\medskip
In particular, one can consider the Pl\"ucker relations that define the Grassmannian, as quadratic relations among minors, coming from
the Laplace expansion of the determinant.

\smallskip
We generalize these classical observations to the tangential and secant variety, by providing analogous local parameterizations. Recall
that the tangential variety is the union of all tangent lines to the variety, while the secant variety is the Zariski closure of the union
of the bisecant lines. It turns out
that the tangential variety is locally isomorphic to a product of an affine space by the Zariski closure of the variety $M$ parameterized by all
minors of degree at least two of a generic matrix. The secant variety is locally isomorphic to a product of an affine space by the cone over $M$.
In particular, the equations of the tangential (resp. secant) variety correspond to (resp. homogeneous) relations among minors of degree at least $2$.

Our method is inspired by the "cumulant trick" coming from statistics. Given probability distributions, the statisticians compute general
moments and cumulants. The formulas for those were the inspiration to define two triangular automorphisms of the
affine space. This method has had other successful applications \cite{ZwiernikSturmfels, MOZ}.

Furthermore, using generalized determinants, we are able to extend our results to all varieties that are both minuscule
and cominuscule obtaining our main theorem. In order to prove our results we present formulas for the generalized determinants for a sum
and a generalized Laplace expansion in Lemmas \ref{lem:detsum} and \ref{lem:laplace}.
The following setting also includes the spinor varieties and the two exceptional Hermitian symmetric spaces.
The equations of the secant and tangential variety correspond to the relations among the Pfaffians.
\begin{tthm}[Theorem \ref{thm:main}]
The secant variety of the minimal orbit $Y$ in the projectivization of a minuscule and cominuscule representation can be covered by
complements of hyperplane sections that are isomorphic to the product of an affine space of the same dimension as $Y$, by the affine cone
over the projective variety parameterized by all generalized determinants of degree strictly greater than one.
\end{tthm}

Analogously, the tangential variety is locally isomorphic to a product of an affine space by the \emph{affine} variety parameterized by all
generalized determinants of degree strictly greater than one. The inclusion of this variety in the cone over it corresponds exactly to the
inclusion of the tangential variety inside the secant.

\begin{cor}
 On each of these affine charts,
 the ideal of the closure of the minimal orbit is given by relations among \emph{all} generalized determinants, while:
\begin{itemize}
\item the secant is defined by the \emph{homogeneous} relations among those generalized determinants that are of degree strictly greater than one,
\item the tangential variety is defined by all the relations among those generalized determinants that are of degree strictly greater than one.
 \end{itemize}
\end{cor}

While the ideal of an equivariantly embedded homogeneous variety is always generated by quadrics, in general,
the ideal of its secant variety is not known. There are no quadrics vanishing on the secant variety. Although, it is not true that the ideal of the secant of a homogeneous variety is always generated by cubics \cite{manivel2009spinor}, it is expected that this property holds for Grassmannians.
We  provide an explicit representation theoretic description of cubics in the ideal of any Grassmannian \ref{thm:cubicinsecant}.
This is related to the plethysm $S^3(\bigwedge^k)$, which has been computed in \cite{littlewood1944invariant, plet2}. We present a short proof in the Appendix.
An analogous description for the secant of Segre varieties was provided in \cite{MR2097214}[Theorem 4.7] and the ideal is now
known in any degree for the secant variety of a Segre-Veronese variety \cite{Raicu}.

There are many motivations to study secant varieties of homogeneous varieties coming both from pure and applied mathematics. Homogeneous
varieties come with a preferred embedding, thus the secant and tangential varieties are intrinsic objects. Their geometry is very interesting - a
classical reference is \cite{Zak}, however the topic has been studied by many other authors - please consult \cite{landsberg2012tensors} and references therein. Still,
the increasing interest in the topic is strongly motivated by possible applications in computer sciences, image processing, statistics etc. This
is related to the problem of determining ranks (or border ranks) of tensors.

The Segre, Veronese embeddings of projective spaces, as well as the Pl\"ucker embeddings of Grassmannians,
are of particular interest as they correspond to general, symmetric, and skew-symmetric tensors.
While the secant and tangential variety of the Segre-Veronse embeddings is well described \cite{Raicu, RaicuOeding}, the problem for Grassmannians is wide open. The description for the secant of $G(3,n)$ was provided in \cite{MR2533306}. In general, we present an easy method to derive the
description of the secant of $G(k,n)$, by the description for $G(k,3k)$ in Proposition \ref{Prop:dimbound}. In \cite{OttavianiLandsberg} one can find a set of equations, so-called skew-flattenings, that define the first and second secant variety on an open subset. 
In \cite{DK2, draisma2013noetherianity} one can find
very nice results concerning bounds for the degree of equations that cut out the secant variety set-theoretically.
The study of relations among minors is also an interesting, difficult topic \cite{bruns2013relations}.

\section{Generalized determinants}

\subsection{Minuscule and cominuscule representations}
Let $\g$ be the Lie algebra of a semisimple complex Lie group $G$. By choosing a Cartan subalgebra $\h$ we obtain a decomposition
$\g=\h\oplus \bigoplus_{\beta\in R} \g_\beta$, where $R$ is the set of roots and the $\g_\beta$ are the root spaces. By fixing
a Weyl chamber, we obtain a basis of the
root system, made of the simple roots $\alpha_1,\dots,\alpha_n$.

Let $\lambda$ be a minuscule and cominuscule dominant weight. In particular $\lambda$ is a fundamental weight $\omega_{i_0}$,
corresponding to the simple root $\alpha_{i_0}$. We get a decomposition $\g=\n\oplus\pp,$ where $\n$ is the sum of the root
spaces $\g_\beta$ where $\beta$ has negative coefficient on $\alpha_{i_0}$ when expressed as a linear combination of the
simple roots. We denote by $B$ the set of these roots. The parabolic subalgebra $\pp$ is the sum of the Cartan subalgebra with the other root spaces.

Let $V_\lambda$ be the minuscule and cominuscule irreducible representation with highest weight $\lambda$.
Its weight decomposition is $V_\lambda=\bigoplus_{w\in F} {W_w}$, where $F$ is the set of weights appearing in $V_\lambda$.
Since $\lambda$ is minuscule $F$ is the orbit of $\lambda$ under the Weyl group action, and  $\dim W_w=1$ for each $w$ in $F$.
By choosing a non zero vector $v_w\in W_w$ we obtain an isomorphism $i_w:W_w\cong \C$. We also define the projection $\pi_w':
V_\lambda\rightarrow W_w$ with respect to the direct sum decomposition,
and the composition $\pi_w=i_w\circ\pi_w'$.

To simplify the notation we define the degree $d(\alpha)=\langle \alpha,\omega_{i_0}^\vee\rangle$, where $\omega_{i_0}^\vee$ is the corresponding
fundamental coweight. In particular, for any $b\in B$ we have $d(b)=1$.
Each element $f\in F$ can be written
as $\lambda$ plus elements of $B$. Such presentation may be not unique, but the number of elements from $B$ is fixed and equals $d(f-\lambda)$.
Moreover, among elements from $\lambda-F$, the maximal degree is $d_{max}=\langle \lambda-w_0(\lambda),\omega_{i_0}^\vee\rangle$, since the lowest weight is the
image of $\lambda$ by the maximal length element $w_0$ of the Weyl group.

There are two infinite series of minuscule and cominuscule representations $V$ of simple complex Lie groups $G$, and two exceptional
examples. They are presented in the table below, with the corresponding values of $d_{max}$. The main geometric
object we are interested in is the minimal orbit $Y$ in $\P (V)$.
$$\begin{array}{ccccl}
 G & V & d_{max} & Y & \\
 SL_n & \wedge^k\c^n & k & G(k,n) & \mathrm{Grassmannian} \\
 Spin_{2n} & \Delta & n & S_{2n} & \mathrm{Spinor\;variety} \\
 E_6 & V_{\omega_1} & 2 & \mathbb{OP}^2 & \mathrm{Cayley\;plane}\\
 E_7 & V_{\omega_7} & 3 & G_\omega(\mathbb{O}^3,\mathbb{O}^6) & \mathrm{Freudenthal\;variety}
\end{array}$$
For Grassmannians we suppose that $k\le n-k$. We will not say much about the two exceptional cases.
The secant variety of the Cayley plane is degenerate, it is the famous Cartan cubic hypersurface
of this Severi variety. The secant of the Freudenthal variety is non degenerate, in fact it is
equal to the whole ambient projective space, and the tangential variety is a quartic hypersurface.

\subsection{Some properties of generalized determinants}
Let us recall the construction of the generalized determinant from \cite{fomin2000recognizing, gel1987combinatorial}. Fix a weight $w\in F$. Define $\det_w:\n\rightarrow\C$ by (cf. diagram \ref{diagram})
$${\det}_w(n)=\pi_w(\exp(n) v_\lambda).$$
\smallskip
Let us make the previous construction more explicit. Let $M_\lambda:=U(\g)\otimes_{U(\b)}\C v_\lambda$ be the Verma module, where $\b\subset \g$ is the
Borel subalgebra and $U$ denotes the universal enveloping algebra. We know by the Poincar\'e-Birkhoff-Witt theorem
that $M_\lambda=\C[(X_b)_{b\in R_-}]v_\lambda$, where $R_-$ is the set of
negative roots. Moreover, $V_\lambda=M_\lambda/K_\lambda$, for $K_\lambda$ the unique maximal submodule. In our case $K_\lambda$ is generated by
relations $X_\beta v_\lambda=0$ for $\beta\in R_-$ such that $\langle\beta,\omega_{i_0}^\vee\rangle =0$, and $X_{\alpha_{i_0}}^2 v_\lambda=0$ \cite{dixmier1977enveloping}.

Notice that so far we have not used the cominuscule assumption. This exactly means that $\n$ is a commutative subalgebra of $\g$. This
allows us to consider the variables $X_\beta$ for $\beta\in B$ as commutative, and treat the $\n$ action on the Verma module just as multiplication of
polynomials.
Once $v_\lambda$ is fixed we may define $v_w$ as follows. For $w=\lambda+\beta$ of degree one,
it is natural to define $v_w$ as the class of $X_\beta v_\lambda$
in $V_\lambda$. For $w$ of higher degree, we \emph{fix} one choice of decomposition
$w=\lambda+\sum_{\beta\in B_w}\beta$, where $B_w\subset B$.
Then we set $$v_w=\Huge(\prod_{\beta\in B_w} X_\beta\Huge) v_\lambda.$$
As for $d(w)\geq 2$ there may be several possible choices,
it is convenient to consider the following compatibility constants.
\dfi[$m(w_1,\dots,w_l)$]
For any  $w_1,\dots,w_l\in F-\{\lambda\}$ such that $\lambda+\sum_{i=1}^l w_i\in F$, we define $m(w_1,\dots,w_l)\in \C$ by the equation
$$m(w_1,\dots,w_l)\Huge(\prod_{i=1}^l\prod_{\beta\in B_{w_i}}X_\beta\Huge) v_\lambda=
\Huge(\prod_{\beta\in B_{\sum_{i=1}^lw_i}}X_\beta\Huge) v_\lambda.$$
\kdfi
\lem[Multiplicative decomposition of compatibility constants]\label{lem:multcompconst}
Consider $\gamma_1,\dots,\gamma_k\in F-\{\lambda\}$ with $\gamma_i=\sum_{j=1}^{a_i} \delta^i_j$, $\delta^i_j\in F-\{\lambda\}$, such that
$\lambda+\sum_{j=1}^k\gamma_j\in F$. Then we have
$$m(\delta^1_1,\dots,\delta^1_{a_1},\delta^2_1,\dots,\delta^k_{a_k})=m(\gamma_1,\dots,\gamma_k)\prod_{i=1}^km(\delta^i_1,\dots,\delta^i_{a_i}).$$
\klem
\dow
Let $\gamma=\gamma_1+\cdots +\gamma_k$. The definition of the compatibility constants implies that
$$\begin{array}{l}
m(\gamma_1,\dots,\gamma_k)(\prod_{i=1}^km(\delta^i_1,\dots,\delta^i_{a_i}))(\prod_{i,j}\prod_{\beta\in B_{\delta^i_j}} X_{\beta})
v_\lambda \\
\hspace{2cm} =  m(\gamma_1,\dots,\gamma_k)(\prod_{i=1}^k\prod_{\beta\in B_{\gamma_i}}X_\beta) v_\lambda =
(\prod_{\beta\in B_{\gamma}}X_\beta) v_\lambda \\
\hspace{5cm} =   m(\delta^1_1,\dots,\delta^1_{a_1},\delta^2_1,\dots,\delta^k_{a_k})(\prod_{i,j}\prod_{\beta\in B_{\delta^i_j}} X_{\beta})v_\lambda,
  \end{array}$$
hence the claim. \kdow

Let us now consider the minimal orbit $Y_\lambda=G[v_\lambda]\subset\P(V_\lambda)$. There is a commutative diagram:
\begin{equation}\label{diagram}\xymatrix{
&&&\P(V_\lambda)&&\\
&G\ar@{->}[r]^{f}&GL(V_\lambda)\ar@{->}[r]^{(v_\lambda)}&V_\lambda\ar@{->}[r]^{\pi'_w}\ar@{-->}[u]&W_w\ar@{=}[r]^{i_w}&\C\\
\n\ar@{^{(}->}[r]&\g\ar@{->}[u]^{exp}\ar@{->}[r]^{df}&gl(V_\lambda)\ar@{->}[u]^{exp}&&\\
}\end{equation}

\medskip
Since $\pp$ is the tangent space to the orbit of $[v_\lambda]$ and $\n$ is transverse to $\pp$, its image in the orbit is a
dense open subset, isomorphic to $\n$ since the latter is nilpotent. Therefore, the generalized determinants $ \pi_w(exp(\cdot))v_\lambda)$
provide a local parametrization of the orbit of the highest weight vector in the given coordinates. Note that
the image of $\n$ is contained in the affine subspace $\a_\lambda\subset \P(V_\lambda)$ given by $x_\lambda\neq 0$. We also see that the
parametrization can be obtained by
$\n\ni n_0\rightarrow \sum_{i=0}^\infty \frac{n_0^i}{i!}v_\lambda\in V_\lambda$, as sufficiently high power of $n_0^i$ acts trivially on
$v_\lambda$, by the weight computation. We have:
\wn[cf. \cite{manivel2009spinor} Sections 2.3, 2.4, 2.5]
The ideal of $\a_\lambda\cap Y_\lambda$ is defined by the relations among all the generalized determinants. In particular,
for the Grassmannian $G(k,n)$ these are relations among all the minors of a generic $k\times(n-k)$ matrix. For spinor varieties
they are the relations between the subpfaffians of a generic skew symmetric matrix.
\kwn
One of the main aims of this article is to generalize the previous theorem to secant varieties. So far we have obtained a parametrization:
$$p':\n\rightarrow \a_\lambda.$$
The coordinates of $\a_\lambda$ are $x_b$ for $b\in F$, $b\neq\lambda$. We have $x_b(p'(n_0))={\det}_b(n_0)$. This provides a
dominant map to the secant variety:
$$p:\C\times \n\times \n\rightarrow \a_\lambda,$$
where $x_b(p((t,n_0,n_1)))=t{\det}_b(n_0)+(1-t){\det}_b(n_1)$.

Let us also present further easy results on generalized determinants. We start by providing a precise formula.
\ob\label{prop:det}
For any weight $\lambda+\gamma\in F$, we have
$${\det}_{\lambda+\gamma}(\sum_{\beta\in B} a_\beta X_\beta)=\sum_{\beta_1+\dots+\beta_{d(\gamma)}=\gamma}
m(\beta_1,\dots,\beta_{d(\gamma)})\prod_{i=1}^{d(\gamma)} a_{\beta_i}.$$
Here we take the sum over decompositions $\beta_1+\dots+\beta_{d(\gamma)}$ considered up to permutations.
In particular, ${\det}_{\gamma+\lambda}$ is a homogeneous polynomial in the $a_\beta$'s, of degree  $d(\gamma)$.
\kob
\dow
We simply expand the exponential of $\sum_{\beta\in B} a_\beta X_\beta$, using the usual formula since
the $X_\beta$ commute. The factor $d(\gamma)!$ appearing in the expansion
simplifies when we consider decompositions of $\gamma$ only up to permutations.
\kdow
By combining Lemma \ref{lem:multcompconst} and Proposition \ref{prop:det} we deduce the following Lemmas.
\lem\label{lem:detsum}
For any weight $\lambda+\gamma\in F$, and any $A,B\in\n$, we have
$${\det}_{\lambda+\gamma}(A+B)=\sum_{\gamma_1+\gamma_2=\gamma}m(\gamma_1,\gamma_2){\det}_{\lambda+\gamma_1}(A){\det}_{\lambda+\gamma_2}(B).$$
\klem
\dow
Apply the previous Lemma, expand the products, use the compatibility conditions, and
finally apply the previous Lemma again.
\kdow

\lem[Generalized Laplace extension]\label{lem:laplace}
Fix a partition $d_1+\dots+d_k=d(\gamma)$, where  $d_1,\dots,d_k\in \Z_+$. Then
$$\binom{d(\gamma)}{d_1,\dots,d_k}{\det}_{\lambda+\gamma}=\sum_{\gamma_1+\dots+\gamma_k=\gamma}m(\gamma_1,\dots,\gamma_k)
\prod_{i=1}^k{\det}_{\lambda+\gamma_i},$$
where the sum is taken over all decompositions such that $d(\gamma_i)=d_i$.
\klem
\dow
By Proposition \ref{prop:det},
$$\sum_{\gamma_1+\dots+\gamma_k=\gamma}m(\gamma_1,\dots,\gamma_k)\prod_{i=1}^k{\det}_{\gamma_i+\lambda}(\sum_{\beta\in B}a_\beta X_\beta)=
\hspace*{5cm}$$
$$\hspace*{2cm}=\sum_{\gamma_1+\dots+\gamma_k=\gamma}m(\gamma_1,\dots,\gamma_k)\prod_{i=1}^k(\sum_{\beta_1^i+\dots+\beta_{d(\gamma_i)}^i=\gamma_i}
m(\beta_1^i,\dots,\beta_{d(\gamma_i)}^i)\prod_{j=1}^{d(\gamma_i)}a_{\beta^i_j}).
$$
By expanding the first product we obtain:
$$
\sum_{\gamma_1+\dots+\gamma_k=\gamma}\sum_{\beta_1^i+\dots+\beta_{d(\gamma_i)}^i=\gamma_i}(m(\gamma_1,\dots,\gamma_k)\prod_{i=1}^k
m(\beta_1^i,\dots,\beta_{d(\gamma_i)}^i))\prod_{j,i}a_{\beta^i_j}.
$$
We can now apply Lemma \ref{lem:multcompconst}. We get:
$$\sum_{\gamma_1+\dots+\gamma_k=\gamma}\sum_{\beta_1^i+\dots+\beta_{d(\gamma_i)}^i=\gamma_i}m(\beta_1^1,\dots,\beta^k_{
d(\gamma_k)})\prod_{j,i}a_{\beta^i_j}.$$
The first two sums give a sum over all decompositions of $\gamma$ into elements from $B$, each decomposition counted $\binom{d
(\gamma)}{d_1,\dots,d_k}$ times. By Proposition \ref{prop:det} this gives the result.
\kdow

\section{The secant and tangent varieties in cumulant coordinates}

\subsection{Minuscule cumulants}\label{sec:mincum}
We introduce two changes of coordinates on $\a_\lambda$.

First, we let $y_{\lambda+\beta}:=x_{\lambda+\beta}$ for $\beta\in B$, and if $d(\gamma)>1$,
$$y_{\lambda+\gamma}:=\sum_{\gamma=\gamma_1+\gamma_2}(-1)^{d(\gamma_2)}m(\gamma_1,\gamma_2)x_{\lambda+\gamma_1}
{\det}_{\lambda+\gamma_2}(\sum_{\beta\in B}x_{\gamma+\beta}X_\beta).$$
Notice that this is an automorphism of $\a_\lambda$ as $x_{\lambda+\gamma}$ appears in $y_{\lambda+\gamma}$ with coefficient $1$ corresponding
to the decomposition $\gamma=\gamma+0$, while all the other $x_b$ that appear satisfy $b>\lambda+\gamma$. Thus the automorphism is triangular.

Second, we let $z_{\lambda+\beta}:=y_{\lambda+\beta}$ for $\beta\in B$,  and for $d(\gamma)>1$,
$$z_{\lambda+\gamma}=\sum_{\gamma_1+\dots+\gamma_k=\gamma}\frac{(-1)^k}{\binom{d(\gamma)}{d(\gamma_1),\dots,
d(\gamma_k)}}m(\gamma_1,\dots,\gamma_k)\prod_{i=1}^ky_{\lambda+\gamma_i}.$$
Here the sum is taken over decompositions such that $d(\gamma_i)\geq 2$ for each $i$. As before we obtain a triangular automorphism.

\subsection{The secant in cumulant coordinates}
We start with the following auxiliary lemma, cf. \cite{MOZ}.

\lem\label{lem:comput}
Let $P_k(t)=(-t)^{k}(1-t)+t(1-t)^{k}$. Then
$$P_k(t)\prod_{i=1}^k(a_i-b_i)=\sum_{A\subset\{1,\dots,k\}}(-1)^{k-|A|}(t\prod_{i\in A}a_i+(1-t)\prod_{i\in A}b_i)\prod_{i\in\{1,\dots,k\}
\setminus A}(ta_i+(1-t)b_i)$$
\klem
\dow
We consider the right hand side (RHS) as a polynomial of degree $k$ in $a_i,b_i$ with coefficients in $k[t]$. By pairing sets $A$ differing only
by given $i_0$ we see that the RHS is zero if $a_{i_0}=b_{i_0}$ for some index $i_0$. Thus we must have a factor $\prod_{i=1}^k(a_i-b_i)$ on the LHS,
and the missing factor is a polynomial $P_k(t)$. The proof
is thus reduced to finding $P_k(t)$ for fixed values of $a_i,b_i$. This can be done by an easy induction - c.f. proof of \cite{MOZ}[Lemma 3.1].
\kdow

We are now ready to present the parametrization of the secant in coordinates $y_b$. Recall we have parametrized it by
mapping $(t,A=\sum a_\beta X_\beta,B=\sum b_\beta X_\beta)\in \C\times \n\times\n$ to
$$(t{\det}_{\gamma+\lambda}(\sum a_\beta X_\beta)+(1-t)
{\det}_{\gamma+\lambda}(\sum b_\beta X_\beta))_{\gamma+\lambda\in F}\in \a_\lambda.$$
\lem\label{lem:seciny}
After the change of coordinates to $y_b$, the parametrization of the secant is given by $y_{\lambda+\beta}=ta_\beta+(1-t)b_\beta$
for $\beta\in B$, and for $d(\gamma)>1$ by:
$$y_{\lambda+\gamma}=P_{d(\gamma)}(t){\det}_{\lambda+\gamma}(\sum_{\beta\in B} (a_\beta-b_\beta)X_\beta).$$
\klem

\dow
For $d(\gamma)>1$,  we compute that $y_{\lambda+\gamma}$ is given by
$$\sum_{\gamma=\gamma_1+\gamma_2}(-1)^{d(\gamma_2)}m(\gamma_1,\gamma_2)(t{\det}_{\lambda+\gamma_1}
(\sum_{\beta\in B} a_\beta X_\beta)+ \hspace*{3cm}$$
$$\hspace*{3cm} +(1-t){\det}_{\lambda+\gamma_1}(\sum_{\beta\in B} b_\beta X_\beta)){\det}_{\lambda+\gamma_2}(\sum_{\beta\in B}(ta_\beta+(1-t)b_\beta)X_\beta).$$
We expand the determinants using Proposition \ref{prop:det}. We get
$$\sum_{\gamma=\gamma_1+\gamma_2}(-1)^{d(\gamma_2)}m(\gamma_1,\gamma_2)(\sum_{\sum_j\beta_j=\gamma_1}m(\beta_1,\dots,
\beta_{d(\gamma_1)})(t\prod_j a_{\beta_j} +(1-t)\prod_j b_{\beta_j})\times$$
$$\times (\sum_{\sum_l{\beta_l'=\gamma_2}}m(\beta_1',\dots,\beta_{d(\gamma_2)})\prod_l(ta_{\beta_l}+(1-t)b_{\beta_l})).$$
Note that the triple sum $\sum_{\gamma=\gamma_1+\gamma_2}\sum_{\sum_j\beta_j=\gamma_1}\sum_{\sum_l{\beta_l'=\gamma_2}}$ is just the sum over
all decompositions $\sum_{\sum_i\beta_i=\gamma}\sum_{A\subset\{1,\dots,d(\gamma)\}}$ by identifying $\gamma_1=\sum_{i\in A}\beta_i$.
Therefore we obtain:
$$\sum_{\sum_i\beta_i=\gamma}\sum_{A\subset\{1,\dots,d(\gamma)\}}(-1)^{d(\gamma)-|A|}
m(\sum_{i\in A}\beta_i,\sum_{i\in\{1,\dots,k\}\setminus A}\beta_i)m((\beta_i)_{i\in A})m((\beta_i)_{i\in \{1,\dots,k\}\setminus A})\times
$$
$$(t\prod_{i\in A} a_{\beta_i} +(1-t)\prod_{i\in A} b_{\beta_i})\prod_{i\in\{1,\dots,k\}\setminus A}(ta_{\beta_i}+(1-t)b_{\beta_i}).$$
By applying Lemma \ref{lem:multcompconst} and Lemma \ref{lem:comput} we get:
$$\sum_{\sum_i\beta_i=\gamma}m(\beta_1,\dots,\beta_{d(\gamma)})P_{d(\gamma)}(t)\prod_i(a_{\beta_i}-b_{\beta_i}),$$
which finishes the proof by Proposition \ref{prop:det}.
\kdow
We now compute the parametrization in the coordinates $z_b$.
\lem\label{lem:secinz}
After the change of coordinates to $z_b$ the parametrization of the secant is given by $z_{\lambda+\beta}=ta_\beta+(1-t)b_\beta$
for $\beta\in B$, and for $d(\gamma)>1$ by
$$z_{\lambda+\gamma}=t(1-t)(1-2t)^{d(\gamma)-2}\det{}_{\lambda+\gamma}(\sum_{\beta\in B} (a_\beta-b_\beta)X_\beta).$$
\klem
\dow
For $d(\gamma)>1$ we have
$$z_{\lambda+\gamma}=\sum_{\gamma_1+\dots+\gamma_k=\gamma}\frac{(-1)^k}{\binom{d(\gamma)}{d(\gamma_1),\dots,
d(\gamma_k)}}m(\gamma_1,\dots,\gamma_k)\prod_{i=1}^k P_{d(\gamma_i)}(t)\det{}_{\gamma_i+\lambda}(\sum_{\beta\in B}
(a_\beta-b_\beta)X_\beta),$$
where the sum is taken over such decompositions that $d(\gamma_i)\geq 2$. This equals:
$$\sum_{d_1+\dots+d_k=d(\gamma)}\frac{(-1)^k}{\binom{d(\gamma)}{d_1,\dots,d_k}}(\prod_{i=1}^kP_{d_i}(t))(
\sum_{\gamma_1+\dots+\gamma_k=\gamma}m(\gamma_1,\dots,\gamma_k)\prod_{i=1}^k\det{}_{\gamma_i+\lambda}(\sum_{\beta\in B} (a_\beta-b_\beta)X_\beta)),$$
where the second sum is taken over such decompositions that $d(\gamma_i)=d_i\geq 2$.
By Lemma \ref{lem:laplace} we obtain:
$$\sum_{d_1+\dots+d_k=d(\gamma)}(-1)^k(\prod_{i=1}^kP_{d_i}(t))\det{}_{\lambda+\gamma}(\sum_{\beta\in B} (a_\beta-b_\beta)X_\beta).$$
Thus the statement is reduced to proving:
$$\sum_{d_1+\dots+d_k=d(\gamma)}(-1)^k(\prod_{i=1}^kP_{d_i}(t))=t(1-t)(1-2t)^{d(\gamma)-2},$$
as $P_1(t)=0$.
We leave it as an exercise, c.f. \cite{MOZ}.
\kdow
As $\det{}_{\lambda+\gamma}$ is a polynomial of degree $d(\gamma)$ we obtain the following result.
\tw\label{thm:main}
The secant variety of the minimal orbit $Y$ in a minuscule and cominuscule representation can be covered by complements of hyperplane sections
isomorphic to a product of an affine space of the same dimension as $Y$, with  the affine cone
over the projective variety parameterized by all generalized determinants of degree strictly greater than one.
\ktw
On such an affine chart, the closure of the minimal orbit is defined by the relations
among \emph{all} the generalized determinants, while the secant is defined by the \emph{homogeneous} relations among
those generalized determinants that are of degree strictly greater than one.

\subsection{The tangent in cumulant coordinates}
Similar algebraic tricks can be applied to a parametrization of the tangential variety of the minimal orbit in a minuscule and cominuscule
representation. However, instead of performing the computations we can rely on the results we have obtained so far.
\lem\label{lem:seclines}
Let us fix two points on the minimal orbit $e^{n_1}v_\lambda,e^{n_1+\epsilon n_2}v_\lambda$ for $n_1,n_2\in \n$ and $\epsilon \in \C$.
The secant line in $\a_\lambda$ joining these points is parameterized by:
$$\C\ni t\rightarrow te^{n_1}v_\lambda+(1-t)e^{n_1+\epsilon n_2}v_\lambda\in \a_\lambda.$$
After the change of coordinates to $z_b$ the parametrization of this line is given by:
$$z_{\lambda+\beta}=t\det{}_{\lambda+\beta}(n_1)+(1-t)\det{}_{\lambda+\beta}(n_1+\epsilon n_2),$$
for $\beta\in B$ and for $d(\gamma)>1$ by:
$$z_{\lambda+\gamma}=t(1-t)(1-2t)^{d(\gamma)-2}\det{}_{\lambda+\gamma}(-\epsilon n_2).$$
\klem
\dow
This is a special case of the computation in Lemma \ref{lem:secinz}.
\kdow
\ob
In the coordinates $z_b$, the tangential variety is a product of the affine space defined by the coordinates $z_{\lambda+\beta}$,
by the affine  variety parameterized for $d(\gamma)>1$ by $z_{\lambda+\gamma}=c\det{}_{\lambda+\gamma}(n)$, with $n\in\n$ and $c$ a fixed constant.
\kob
\dow
Let $T$ be the Zariski closure of the given parametrization. First we show that the tangential variety contains $T$. By Proposition \ref{prop:det}
the generalized determinants are homogeneous polynomials. Thus, the secant lines from Lemma \ref{lem:seclines} can be parameterized for
$d(\gamma)>1$ by
$$z_{\lambda+\gamma}=(-\epsilon t)(-\epsilon+\epsilon t)(-\epsilon+2\epsilon t)^{d(\gamma)-2}\det{}_{\lambda+\gamma}(n_2).$$
Taking $t=1/\epsilon$, as $\epsilon\rightarrow 0$ we obtain:
$$z_{\lambda+\gamma}=c\det{}_{\lambda+\gamma}(-2n_2).$$
Notice that the coordinates $z_{\lambda+\beta}$ are arbitrary as $n_1$ can be chosen arbitrary and the generalized determinants of degree one
$\det_{\lambda+\beta}(n_1)$ just return a coordinate of $n_1$. As $\epsilon\rightarrow 0$ the secant lines approach tangent lines, so $T$ is
indeed contained in the tangential variety. Notice that by Theorem \ref{thm:main} the secant variety is a cone over $T$, thus $T$ is of codimension
one in the secant variety, as is the tangential variety. The proposition follows. \kdow

\section{Secants of Grassmannians}
The most classical minimal orbits of representations that are minuscule and cominuscule are Grassmannians. The equations of the secant of the
Grassmannian $G(k,V)$ embedded by the Pl\"ucker embedding are not known in general, apart from the cases $k=2,3$ \cite{MR2533306}. Let us start with
the case $k=2$.
\subsection{Pfaffian-Pl\"ucker and the $\sigma(G(2,n))\sim G(2,n-2)$ correspondence}
Let $V$ be an $n$ dimensional vector space. Consider $G(2,V)\subset \p(\bigwedge^2 V)$. We may represent elements of $\bigwedge^2 V$ as skew
symmetric matrices. The points of $G(2,V)$ correspond to matrices of rank $2$, hence to degree two minors of a
generic $2\times n$ matrix. Locally, $G(2,V)$ can be parameterized by all minors of a generic $2\times (n-2)$ matrix.
The points of the secant $\sigma(G(2,V))$ correspond to matrices of rank at most $4$.
Moreover, the ideal of the secant variety is generated by $6\times 6$
subpfaffians of the generic skew-symmetric matrix.

The Main Theorem \ref{thm:main}
asserts that $\sigma(G(2,n))$ is covered by affine open subsets which are products of a $2(n-2)$-dimensional affine space by the affine cone over
the projective variety $Y$ parameterized by $2\times 2$ minors of a generic $2\times (n-2)$ matrix, which is $G(2,n-2)$. The unique singular
point over the affine cone corresponds exactly to the points in $G(2,n)$ itself. In particular we see that the singularities of $\sigma(G(2,n))$ are
exactly the same as those of the affine cone over $G(2,n-2)$, which have been well-studied. 
\wn
The secant $\sigma(G(2,n))$ is covered by affine varieties $\a^{2(n-2)}\times \widehat G(2,n-2)$, where $\widehat G(2,n-2)$ is the affine cone over the Grassmannian $G(2,n-2)$. In particular it has the same singularities as the cone, e.g. rational -- cf. \cite{MR2533306}.
\kwn
Let us further investigate the correspondence. Fix a basis $e_1,\dots,e_n$ of $V$. In order to fix a principal affine open subset of
$\p(\bigwedge^2 V)$ we have to choose two vectors $e_{i_1},e_{i_2}$. The equations of $G(2,n-2)$ are known as Pl\"ucker relations. They are
indexed by the choice of four further vectors $e_{i_3},e_{i_4}, e_{i_5}, e_{i_6}$. The chosen Pl\"ucker relation after the isomorphisms
described in Section \ref{sec:mincum} has to induce an equation of $\sigma(G(2,n))$. Pl\"ucker relations are quadratic, and so is our change of
coordinates. 
So we could expect a degree $4$ equation. However, it turns out that we obtain the subpfaffian indexed by $i_1,\dots,i_6$,
which has degree three!
\ex
Without loss of generality we assume that $i_1=1$ and $i_2=2$, that is we dehomogenize with respect to the variable $x_{12}$.
The Pl\"ucker relation corresponding to $i_3=3,i_4=4,i_5=5,i_6=6$ is given by:
$$z_{34}z_{56}+z_{45}z_{36}-z_{46}z_{35}.$$
The isomorphism of affine spaces is defied in this case by:

$z_{ij}=x_{ij}$ for $i\le 2<j$, and $z_{ij}=x_{12}x_{ij}-x_{1i}x_{2j}
+x_{1j}x_{2i}$ for $2<i<j$.

After the substitution the Pl\"ucker relation equals the Pfaffian of the matrix:
$$ \left( \begin{array}{cccccc}
0 & x_{12} & x_{13} & x_{14} & x_{15} & x_{16} \\
-x_{12} & 0 & x_{23} & x_{24} & x_{25} & x_{26} \\
-x_{13} & -x_{23} & 0 & x_{34} & x_{35} & x_{36} \\
-x_{14} & -x_{24} & -x_{34} & 0 & x_{45} & x_{46} \\
-x_{15} & -x_{25} & -x_{35} & -x_{45} & 0 & x_{56} \\
-x_{16} & -x_{26} & -x_{36} & -x_{46} & -x_{56} & 0 \\
\end{array} \right).$$
\kex

\subsection{General case}
Motivated by Main Theorem \ref{thm:main} we make the following definition:
\dfi[$\tilde M_{a,b}$, $M_{a,b}$]
For two natural numbers $a,b\in \Z$ we define the embedded projective variety $\tilde M_{a,b}$ as the Zariski closure of the map that to a
generic $a\times b$ matrix associates all its minors of degree at least $2$.
In particular, $\tilde M_{2,b}=G(2,b)$. Let $M_{a,b}$ be the affine cone over $\tilde M_{a,b}$.
\kdfi
As another corollary of \ref{thm:main} we obtain:
\wn
The secant variety $\sigma(G(k,n))$ is covered by principal open affine sets isomorphic to $\a^{k\times(n-k)}\times M_{k,n-k}$.
\kwn
Thus the equations of $\sigma(G(k,n))$, after the affine isomorphism, give homogeneous relations among minors and vice versa. The study of
relations among minors is a very classical problem. In general the generators of these relations are unknown, apart from the case of maximal
minors, when we have Pl\"ucker relations. Still, it is very easy to produce some quadratic relations by Laplace expansion. Moreover, the
conjectural generators in case of $2\times 2$ minors were provided in \cite{bruns2013relations}.

\section{The secant and representation theory}

In this section we consider the equations of secant varieties of flag manifolds from the perspective
of representation theory. The main case we have in mind is that of Grassmannians.

Let $\lambda$ be a Young diagram with $k$ rows. Let $Y_\lambda$ be the unique closed $GL(V)$-orbit in $\p(S_\lambda V)$,
which is the orbit of highest weight vectors.

\subsection{Equations of the secant variety}\label{subsec:Eqdefsec}
We start by recalling a method for describing the equations of secant varieties due to Landsberg and Manivel \cite{MR2097214}.

A polynomial $P$ belongs to the ideal of the $(s-1)$-st secant variety if and only if for any points $Q_1,\dots, Q_s$ in the affine cone over $Y_\lambda$,
which we call simple vectors,
and any complex numbers $a_1,\dots,a_s$ we have $P(a_1Q_1+\dots+a_sQ_s)=0$. Suppose that $P$ is homogeneous of degree $d$. Recall that to any homogeneous
polynomial $P$ of degree $d$ one can associate a $d$-linear form $\widetilde P$ called the polarization of $P$, 
so that $\widetilde P(v,\dots,v)=P(v)$. It can be defined as:
$$\widetilde P(v_1,\dots,v_d)=\frac{1}{d!}\frac{\partial}{\partial \lambda_1}\dots\frac{\partial}{\partial \lambda_d}
P(\lambda_1v_1+\dots+\lambda_dv_d)_{|\lambda=0}.$$
If we consider $P(a_1Q_1+\dots+a_sQ_s)$ as a polynomial in the $a_i$'s, the coefficient of the monomial $\prod a_i^{\alpha_i}$, up to
some constant, is equal to $\widetilde P(Q_1,Q_1,\dots,Q_s)$, where each $Q_i$ appears $\alpha_i$ times. 
All these polynomials must therefore vanish.

The vanishing of a polarization on all $d$-tuples of simple vectors implies its vanishing on the whole space (since $Y_\lambda$ is
linearly non degenerate), so this immediately implies:
\wn
There are no non zero polynomials of degree $d$ vanishing on the $(s-1)$-st secant variety for $d\leq s$.\kwadrat
\kwn
\dfi[$l\cdot\lambda$]
Let $l\cdot\lambda$ be the Young diagram, obtained by multiplying each row of $\lambda$ by $l$. The corresponding representation
$S_{l\cdot\lambda}V$ is spanned by the (cone over the) $l$-th Veronese embedding of $Y_\lambda$.
\kdfi
Let us consider a decomposition $\sum_{i=1}^s \alpha_i=d$. To this decomposition we associate the representation $M$ defined by
$$M:=S_{\alpha_1\cdot\lambda}V\otimes\dots\otimes S_{\alpha_s\cdot\lambda}V.$$
The vector space $M$ is spanned by simple tensors of the form
$v_{\alpha_1}(Q_1)\otimes\dots\otimes v_{\alpha_s}(Q_s)$, where $v_i$ is the $i$-th Veronese embedding and the $Q_i$ are simple vectors
in $S_\lambda V$. We may define a linear form on
$S^d(S_\lambda V^*)\otimes M$,
that to a polynomial $P$ of degree $d$ and a tensor $v_{\alpha_1}(Q_1)\otimes\dots\otimes v_{\alpha_s}(Q_s)\in M$ associates $\widetilde
P(Q_1,Q_1,\dots,Q_s)$, where each $Q_i$ appears $\alpha_i$ times. After dualizing,
we get a map $$f_{\alpha} :S^d(S_\lambda V^*)\rightarrow M^*.$$
The observations above can be stated as follows.

\ob\label{kernel}
Let $I_d$ be the $d$-th graded part of the ideal of $\sigma_{s-1}(Y_\lambda)$. Then
$$I_d=\cap_\alpha \ker f_\alpha,$$
where the intersection is over all possible decompositions $\sum_{i=1}^s \alpha_i=d$.\kwadrat
\kob
The kernels of the maps $f_{\alpha}$ are hard to describe in general.
However, we can restrict each $f_\alpha$ to an isotypic component $R_\mu$ corresponding to a Young diagram $\mu$ to obtain a bound on its multiplicity
in the ideal.
\wn\label{crucial}
Let us fix a component $\mu$ of weight $dk$. For a given decomposition $\alpha$ given by $\sum_{i=1}^s \alpha_i=d$ let $m_{\alpha}$ be the
multiplicity of the component given by $\mu$ inside $\bigotimes_{i=1}^s S_{\alpha_i\cdot\lambda}V^*$. Let $m=\sum_{\alpha} m_{\alpha}$, where the
sum is over all possible decompositions $\alpha$ of $d$.

Then the multiplicity of $S_\mu V^*$ in the degree $d$ part of
the ideal of $\sigma_{s-1}(Y_\lambda)$ is at least its multiplicity inside $S^d(S_\lambda V^*)$ minus $m$.\kwadrat
\kwn
Let us give a first application. Recall that $\lambda$ has $k$ rows.
\wn\label{cor:longarein}
All isotypic components of $S^d(S_\lambda V^*)$ that correspond to Young diagrams with either
\begin{itemize}
\item the $\ulcorner \frac{d}{s}\urcorner$ column of length strictly less than $k$,
\item more than $ks$ rows,
\end{itemize}
are contained in the ideal of $\sigma_{s-1}(Y_\lambda)$.
\kwn

\dow
Under each of these two hypotheses, it follows from the Littlewood-Richardson rule that $m_\alpha=0$ for
any decomposition $\alpha$ of $d$.
\kdow

\subsection{Stabilization}
To a vector space $V$ and a Young diagram $\lambda$ with $k$ rows we have associated a projective variety $Y_\lambda$. In this subsection
the vector space $V$ will not be fixed, thus we will consider the variety $Y_{\lambda,V}$. As the dimension of $V$ changes, also the properties
of the secant variety $\sigma_s(Y_{\lambda,V})$ may change. However, certain properties will stabilize, as the dimension of $V$ grows.
\ob
The multiplicities of isotypic components of the graded ring of $\sigma_s(Y_{\lambda,V})$, considered as $GL(V)$-representations,
are independent of $V$ as long as $\dim V\geq (s+1)k$.
\kob

\dow
This follows from the methods explained in Subsection \ref{subsec:Eqdefsec}. Indeed, in the domains of the maps $f_{\alpha}$,
there only appear representations given by Young diagrams with at most $(s+1)k$ rows.
\kdow
Among interesting properties of $\sigma_s(Y_\lambda)$ such as being (arithmetically) Gorenstein, (arithmetically) Cohen-Macaulay,
(projectively) normal etc., we do not know which   ones stabilize and under which bounds on the dimension. As far as equations
are concerned we have the following result.

\ob\label{Prop:dimbound}
Suppose that  $\sigma_s(Y_{\lambda,V})$ is defined ideal-theoretically (resp. set-theoretically, resp. scheme-theoretically) by
equations of degree at most $d$ for $\dim V=(s+2)k$. Then the same is true for $\sigma_s(Y_{\lambda, V'})$ for any $V'$
with $\dim V'\geq\dim V$.
\kob
\dow
We prove only the ideal-theoretic version, the other versions being analogous.
Let $I_V$ and $I_{V'}$ be the ideal of respectively $\sigma_s(Y_{\lambda,V})$ and
$\sigma_s(Y_{\lambda,V'})$, where $\dim V'\geq\dim V$. We consider $V$ as a subspace of $V'$ and we choose
a splitting $V'=V\oplus W$. Then there is an induced projection map from $S_\lambda V'$ to  $S_\lambda V$,
and the image $Y_{\lambda,V'}$ is exactly $Y_{\lambda,V}$. Therefore the same is true for the secants
and this implies that  $I_V$ is contained in $I_{V'}$.

We choose a maximal torus in $GL(V')$ given by the product of a maximal torus in $GL(V)$ times a maximal torus in $GL(W)$.
Moreover we order a basis in which this torus is made of diagonal matrices by taking vectors from $V$ before wectors from $W$.
This implies that a highest weight vector in any $S_\mu V'$ belongs to $S_\mu V$ if  $\mu$ has less than $d$ rows.

We will use this observation to show that all highest weight vectors in $I_{V'}$ are generated in degree at most $d$.
As the ideal is $GL(V')$-equivariant the  proposition will follow.

\smallbreak\noindent {\it Step 1}. We first consider highest weight vectors corresponding to Young diagrams with at most $d=k(s+2)$ rows.
As we have just seen, a polynomial $P\in I_{V'}$ that is such a highest weight vector belongs to $S^*(S_\lambda V)$.
Therefore $P$ belongs to $I_V$, thus it is generated in degree
at most $d$ in $I_V$ over the ring $S^*(S_\lambda V)^*$, hence also in $I_{V'}$ over the ring $S^*(S_\lambda V')^*$ since $I_V\subset I_{V'}$.


\smallbreak\noindent {\it Step 2}. Now we consider highest weight vectors corresponding to Young diagrams with more than $k(s+2)$ rows.
Fix a polynomial $P\in I_{V'}$, of degree $d'$, that is such a highest weight vector. The proof is inductive on $d'$, which we
may suppose to be bigger than $d$.
Consider the canonical map:
$$mult : S^{d'-1}(S_\lambda V^*)\otimes S_\lambda V^*\rightarrow S^{d'}(S_\lambda V^*),$$
that corresponds to multiplication of polynomials.
We may decompose the representation $S^{d'-1}(S_\lambda V^*)$ into a direct sum of isotypic components. Note that by the Littlewood-Richardson
rule, only isotypic components with more than $k(s+1)$ rows can be mapped by $mult$ to the isotypic component represented by $\mu$. As $mult$
is surjective we see that:
$$P=\sum l_iQ_i,$$
where the $l_i$ are linear forms and the $Q_i$ belong to isotypic components represented by Young diagrams with more than $k(s+1)$ rows. By Corollary
\ref{cor:longarein} we know that $Q_i$ belongs to $I_{V'}$. By the inductive assumption all $Q_i$ are generated in degree at most $d$, which finishes the proof.
\kdow

\subsection{Cubic equations of the secant}
In this subsection we will provide the description of cubics vanishing on the secant of a Grassmannian $G(k,V)$
as a subrepresentation of $S^3(\bigwedge^k V)^*$.
The space of such cubics will be denoted by $I_3$. We will give explicit formulas for the multiplicities of each irreducible component.
When these multiplicities are smaller than the corresponding multiplicities in $S^3(\bigwedge^k V)^*$ we will explicitly provide linear forms
cutting out the highest weight space in $I_3$. A partial result was given in \cite[Proposition 1.6, Section 3]{MR2533306}.

Let us start with general remarks. First note that Lemma \ref{redukcja} applies to the ideal of the secant as follows.
\ob[Reduction for secants]\label{redsecgras}
Let $\lambda$ be a Young diagram with $d$ columns, and denote by $\lambda'$ the diagram obtained by removing the first row.
Then the multiplicity of the component corresponding to $\lambda$ inside $I_d(\sigma_s(G(k,n)))$ is equal to the
multiplicity of the component corresponding to $\lambda'$  inside $I_d(\sigma_s(G(k-1,n-1)))$.
\kob
\dow
Let $P\in S^d(\bigwedge^k \c^n)^*$ be any polynomial in the highest weight space corresponding to $\lambda$. Denoting by $e_1,\ldots ,e_n$ the
canonical basis of $\c^n$, let $Q$ be the set of variables that
correspond to wedge products of basis elements that contain $e_1$. The polynomial $P$ is a polynomial in these variables only. Let $i$ be an
application that to a wedge product of basis vectors containing $e_1$ associates the same wedge product without $e_1$ and with all indices
decreased by one. We can identify the variables from $Q$ with variables of $\bigwedge^{k-1}\c^{n-1}$. By this
identification points of the Grassmannian $G(k,n)$ correspond to points of $G(k-1,n-1)$. Moreover
$P$ belongs to the ideal of $\sigma_s G(k,n)$ if and only if $i(P)$ belongs to the ideal of $\sigma_s G(k-1,n-1)$.
\kdow

Let us focus on the case where $s=2$ and $d=3$.
Using Proposition \ref{redsecgras} it is enough to obtain the multiplicities corresponding to Young diagrams with at most two columns.
For this we apply Proposition \ref{kernel}.
There are only two decompositions of $3$ to consider: $3=1+2$ and $3=0+3$. The latter provides only one component, corresponding
to the Young diagram with three columns of length $k$. The decomposition $3=1+2$ gives many components, but only one of them has two columns.
This component has got the first column of length $2k$ and the second one of length $k$. It appears with multiplicity one in the
tensor product.

Using Lemma \ref{crucial} we get the following Proposition.
\ob
All representations corresponding to Young diagrams with two columns different from $(2k,k)$ appear with the same multiplicity in
$S^3(\bigwedge ^k V^*)$ as in $I_3$. For the component $(2k,k)$ the multiplicity may drop at most by $1$.
\kwadrat
\kob
We will now prove that in fact the latter multiplicity always drops by one. We will do this by providing a linear form on the highest weight space
that will cut out the ideal. More precisely we will give an example of a polynomial in the highest weight space of $(2k,k)$ that does not vanish on the point
$Q:=e_1\wedge\dots\wedge e_k+e_{k+1}\wedge\dots\wedge e_{2k}$. Our construction is motivated by a general method of defining highest weight
vectors in plethysms described in \cite{JaLaurent}.
\dfi[$x_{a_1\dots a_k}$]
By $x_{a_1\dots a_k}$ we denote a linear form, or equivalently a variable, corresponding to $e_{a_1}^*\wedge\dots\wedge e_{a_k}^*$. If we
permute the $a_i$'s then the variable changes sign according to the sign of the permutation.
\kdfi
For $k$ even we consider the polynomial $$P=\sum_{\sigma\in S_{2k}}\sgn(\sigma) x_{\sigma(1)\dots\sigma(k)}x_{\sigma(k+1)\dots\sigma(2k)}x_{1\dots k}.$$
This is a highest weight vector for the weight $(2k,k)$.
Moreover it contains only one monomial that is nonzero on $Q$, namely $x_{1\dots k}^2 x_{k+1\dots 2k}$. In particular $P(Q)\ne 0$.

For $k$ odd we consider the polynomial
$$P=\sum_{\sigma\in S_{2k}, \delta\in S_k}\sgn(\sigma) \sgn(\delta) x_{\sigma(1)\dots \sigma(k)} x_{\sigma(k+1) \dots\sigma(2k-1),\delta(1)}
x_{\sigma(2k),\delta(2)\dots\delta(k)}.$$
We want to show that $P(Q)\neq 0$. First let us consider one monomial $$x_{\sigma(1)\dots \sigma(k)} x_{\sigma(k+1) \dots\sigma(2k-1),\delta(1)}
x_{\sigma(2k),\delta(2)\dots\delta(k)}$$ appearing in the sum defining $P$. If the monomial is nonzero on $Q$ all the variables, up to
permutation of indices must be either $x_{1\dots k}$ or $x_{k+1\dots 2k}$. As $x_{\sigma(k+1) \dots\sigma(2k-1),\delta(1)}$ and $x_{\sigma(2k),
\delta(2)\dots\delta(k)}$ contain indices less or equal to $k$ they must be, up to sign, $x_{1\dots k}$. Hence $x_{\sigma(1)\dots \sigma(k)}$
must be $x_{k+1\dots 2k}$. We see that $P(Q)\neq 0$ is equivalent to $$\sum_{\sigma\in S_{k}, \delta\in S_k,\sigma(k)=\delta(1)}\sgn(\sigma)
\sgn(\delta) x_{\sigma(1) \dots\sigma(k-1),\delta(1)}x_{\sigma(k),\delta(2)\dots\delta(k)}x_{k+1,\dots, 2k}\neq 0.$$ Notice that since $k$ is odd
all these monomials have coefficients of the same sign, hence  the sum is nonzero.
This concludes the proof of  the main theorem of this subsection:
\tw\label{thm:cubicinsecant}
Let $(a,b,c)$ denote the isotypic component corresponding to the Young diagram with three columns of lengths respectively $a,b,c$.
The multiplicity of the component $(a,b,c)$ in the ideal $I_3(\sigma(G(k,n)))$ is zero for $n<a$. Otherwise it is
equal to
\begin{enumerate}
\item the multiplicity of $(a-c,b-c,0)$ in $S^3(\bigwedge^{k-c})$ if $a-c\neq 2(b-c)$,
\item the multiplicity of $(2(b-c),b-c,0)$ in $S^3(\bigwedge^{k-c})$ minus one if $a-c=2(b-c)$ (and $b=k$). In this case the polynomial in the
highest weight space is in $I_3$ if and only if it does not contain the monomial $x_{1,\dots,k}^2x_{k+1,\dots,2k}$.
\end{enumerate}
\ktw
In other words, by restricting cubics to the secant variety $\sigma(G(k,n))$ we get
$$\c [\sigma(G(k,n))]_3=\bigoplus_{c=0}^k S_{\alpha(k,c)}V^*,$$
where the Young diagram $\alpha(k,c)$ has colums of lengths $(2k-c,k,c)$.

\subsection{The complexity of the secant}

Recall that the complexity of a $G$-variety $X$, where $G$ is a reductive group, is defined as
the codimension of the generic $B$-orbit, where $B$ is a Borel subgroup of $G$.
By \cite{MR2363430} the complexity of the tangential variety is zero (otherwise said the tangential
variety is spherical). This immediately implies that the secant variety has complexity
at most one. Let us prove that there is equality.

\ob
The complexity of the secant variety of a Grassmannian is one.
\kob

\dow
For simplicity we just treat the case of a Grassmannian $G(k,2k)$, $k\ge 3$, the general case being similar.
A generic point of the secant is of the form $p=[u_1\wedge\cdots\wedge u_k+v_1\wedge\cdots\wedge v_k]$,
where $u_1,\ldots ,u_k$ and $v_1,\ldots ,v_k$ are basis of two transverse subspaces $U$ and $V$.
Suppose that $M\in GL(2k)$ belongs to the connected component of the stabilizer of $p$, then we have
$M=X+Y$ where $X\in End(U)$ and $Y\in End(V)$ are such that $\det(X)=\det(Y)$.

Consider a Borel subgroup $B$ of $GL(2k)$ defined as the stabilizer of a generic complete
flag of subspaces $L_1\subset\cdots\subset L_{2k-1}$. Each $L_i$ is generated by vectors $\ell_1,
\ldots ,\ell_i$ and we may suppose that for $i\le k$, $\ell_i=u_i+v_i$. For $i>k$, we may suppose
that $\ell_i=a_{i-k}$ belongs to $U$ or $\ell_i=b_{i-k}$ belongs to $V$.

If $M\in B$ belongs to the connected component of the stabilizer of $p$, we decompose $M=X+Y$ as above.
Then $X$ (respectively $Y$) must preserve the flag of subspaces of $U$ (respectively $V$) defined
by the vectors $u_1,\ldots ,u_k$ (respectively $v_1,\ldots ,v_k$), and the matrices of $X$ and $Y$
in these basis must be the same. Moreover, $M$ also has to preserve the flag defined by the vectors
$a_1,\ldots ,a_k$. But the intersection in $GL(k)$ of two Borel subgroups stabilizing two complete
flags in general position is just a maximal torus. Therefore the connected component of the stabilizer
of $p$ in $B$ is isomorphic to a maximal torus of $GL(k)$, in particular its dimension is $k$.

We conclude that the $B$-orbit of $p$ has dimension $\dim(B)-k=2k^2$ which is one less than the
dimension of the secant variety.
\kdow

A consequence of this observation is that the multiplicities in the coordinate ring of the secant
variety can only grow linearly. These multiplicities are bounded by those of the coordinate ring
of the $GL(2k)$ orbit $\mathcal{O}$ of $p$, which is open in the secant. We have seen in the previous proof that
the connected component of the stabilizer is the subgroup $S(GL(k)\times GL(k))$ of $GL(k)\times GL(k)$
defined as the set of pairs of matrices with the same determinant. More precisely
$$\mathcal{O}\simeq GL(2k)/S(GL(k)\times GL(k)) \rtimes \mathbb{Z}_2.$$

\ob\label{prop:decomporbit}
The coordinate ring of the open orbit $\mathcal{O}$ in the secant variety of the Grassmannian $G(k,W)$,
where $W$ has dimension $2k$, is
$$\mathbb{C}[\mathcal{O}] = \bigoplus_{\alpha\in D_k}\lceil\frac{\alpha_k-\alpha_{k+1}+1}{2}\rceil S_\alpha W^*,$$
where $D_k$ is the set of non increasing sequences $\alpha=(\alpha_1,\ldots ,\alpha_{2k})$ of relative
integers, such that $\alpha_i+\alpha_{2k+1-i}$ is independent of $i$. Here, $\alpha_i$ is the length of the $i$-th \emph{row} of the corresponding Young diagram.
\kob

\dow
With the same notation as before, the $GL(2k)$-orbit of the point $u_1\wedge\cdots\wedge u_k+v_1\wedge\cdots\wedge v_k$
is open in the cone over $\mathcal{O}$. Its stabilizer is $L=(SL(k)\times SL(k)) \rtimes \mathbb{Z}_2$, with
connected component  $L^0=SL(k)\times SL(k)$ . By the Peter-Weyl theorem we have
$$\mathbb{C}[\mathcal{O}] = \mathbb{C}[G]^{L} =\bigoplus_{\alpha}\dim (S_\alpha W)^L \; S_\alpha W^*,$$
where the sum is over all the non increasing sequences $\alpha=(\alpha_1,\ldots ,\alpha_{2k})$ of relative
integers. In order to determine the dimension of the space $(S_\alpha W)^L$ of $L$-invariants
we first consider the $L^0$-invariants. Let us write $\alpha_i=\lambda_i+\ell$, where $\lambda$ is a partition
with $2k$-th part equal to zero, and $\ell\in\mathbb{Z}$. In the decomposition formula
$$S_\alpha W = S_\lambda W\otimes (\det W)^\ell=\bigoplus_{\mu,\nu}c^\lambda_{\mu,\nu}S_\mu U\otimes S_\nu V\otimes
(\det U)^\ell\otimes (\det V)^\ell,$$
where the $c^\lambda_{\mu,\nu}$ are the Littlewood-Richardson coefficients, we see that in order to get
$L^0$-invariants we need to take $\mu=(m^k)$ and $\nu=(n^k)$ for some integers $m$ and $n$ (where by $(m^k)$
we mean the partition with $k$ parts equal to $m$). Then $S_\mu U\otimes S_\nu V= (\det U)^m\otimes (\det V)^n$,
and we will get a one-dimensional space of $L^0$-invariants.

The Littlewood-Richardson rule shows that for $c^\lambda_{\mu,\nu}$ to be non zero, the partition $\lambda$ must
be of form $(m+\theta_1,\ldots,m+\theta_k,n-\theta_k,\ldots ,n-\theta_1)$. In particular, since $\lambda_{2k}=0$,
$\theta_1=n$. Moreover $\lambda_i+\lambda_{2k+1-i}=m+n$ is independent of $i$. If these conditions are fulfilled,
then $c^\lambda_{\mu,\nu}=1$. Note that $\lambda$ being given, there are several possibilities for $m=n$,
subject to the constraints that $m+n=\lambda_1$, $m\ge\lambda_{k+1}$ and $n\ge 0$. This means that $\lambda_{k+1}\le
m\le\lambda_1$, and $0\le n\le\lambda_k$.

Now we can deduce the $L$-invariants. Indeed the $\mathbb{Z}_2$ factors in $L$ switches $U$ and $V$, hence $m$ and $n$.
In particular $m$ and $n$ will only contribute to the $L$-invariants in the range $\lambda_{k+1}\le m,n \le\lambda_k$,
and by symmetry $(m,n)$ and $(m',n')=(n,m)$ contribute to a single $L$-invariant. This implies the statement.
\kdow
\uwa
Let us note that the formula for the multiplicities of the isotypic components in $\mathbb{C}[\mathcal{O}]$ from Proposition \ref{prop:decomporbit} exactly coincides with the upper bounds for the multiplicities in the algebra of the secant variety obtained in Corollary \ref{crucial}. Indeed, for those isotypic components each $m_\alpha$ in the Corollary equals one.
\kuwa

\section{Appendix - Plethysm}\label{sec:app2}

We are investigating the space of polynomials vanishing on the secant of a Grassmannian $G(k,V)$.
Thus the representation $S^d(\bigwedge^k V)^*$  is of great importance for us.
Its decomposition is not known in general. However, it is know for $d\leq 3$ \cite{plet2}. For the sake of
completeness, and as the results we found contained some misprints we present an easy, combinatorial proof.

We will be using the following duality result.
\fa[\cite{Carre}, \cite{MR1651092}]\label{dualities}
$$S^{\mu}(S^{2l} V)=S^{\mu}(\bigwedge^{2l} V)^\vee,\qquad S^{\mu}(S^{2l+1} V)=S^{\mu^\vee}(\bigwedge^{2l+1}V)^\vee,$$
where $^\vee$ means that each irreducible component corresponding to a Young diagram $\nu$ is replaced with the component corresponding to the
transpose of $\nu$, denoted $\nu^\vee$.
\kfa

\tw\label{th:plethysm3}
The multiplicity in $S^3(\bigwedge^k)$ of the isotypic component corresponding to the Young diagram with columns $(a,b,c)$, where $a+b+c=3k$, equals:
\begin{enumerate}
\item for $\min(b-c,a-b)$ is even:

if $\max(b-c,a-b)$ is even then $\ulcorner\frac{\min (b-c+1, a-b+1) }{6}\urcorner$,

if $\max(b-c,a-b)$ is odd then $\llcorner\frac{\min (b-c+1, a-b+1) }{6}\lrcorner$,
\item for $\min(b-c,a-b)$ is odd

if $\min(b-c,a-b)=0$ mod $3$ then $\ulcorner\frac{\min (b-c+1, a-b+1) }{6}\urcorner$,

if $\min(b-c,a-b)=1$ mod $3$ then $\llcorner\frac{\min (b-c+1, a-b+1) }{6}\lrcorner$,

if $\min(b-c,a-b)=2$ mod $3$ then $\frac{\min (b-c+1, a-b+1) }{6}$.
\kwadrat
\end{enumerate}

\ktw
The following Lemma \ref{redukcja} and Corollary \ref{Redlemsym} are classical. A variation of them can be found for example in
\cite[5.8, 5.9]{Carre}. However, the proofs that we know usually take advantage of properties of Schur polynomials. We propose a
very simple, direct approach, that not only provides equality of multiplicities of isotypic components, but also explicitly gives an isomorphism.
\lem[Reduction Lemma]\label{redukcja}
Let $\mu$ be any Young diagram of weight $n$. Let $\lambda$ be a Young diagram with $n$ columns and weight $nk$. Let $\lambda'$ be $\lambda$
with the first row removed. The multiplicity of the component corresponding to $\lambda$ in $S^{\mu}(\bigwedge^k W)$ equals the multiplicity
of the component corresponding to $\lambda'$ in $S^{\mu}(\bigwedge^{k-1} W)$.
\klem
\dow
Consider the inclusion $S^{\mu}(\bigwedge^k V)\subset (\bigwedge^k V)^{\otimes n}$ with a basis given by tensor products of wedge product of
basis elements of $V$. Each vector in the highest weight space corresponding to $\lambda$ must contain exactly one $e_1$ in each tensor. We
get an isomorphism of highest weight spaces by removing $e_1$ and decreasing by one the indices of other basis vectors.
\kdow
Due to the dualities \ref{dualities} we get the following corollary.
\wn\label{Redlemsym}
Let $\mu$ be any Young diagram of weight $n$. Let $\lambda$ be a Young diagram with $n$ rows and weight $nk$. Let $\lambda'$ be equal to
$\lambda$ with the first column removed. The multiplicity of the component corresponding to $\lambda$ in $S^{\mu}(S^k W)$ equals the
multiplicity of the component corresponding to $\lambda'$ in $S^{\mu^\vee}(S^{k-1} W)$.
\kwn
Let us give some applications of these easy observations.
First we prove the classical decompositions, first obtained by Thrall \cite{plet1}, \cite{plet2}[4.1-4.6]:
\ob\label{drugasym}
One has $Gl(W)$-modules decompositions
$$S^2(S^n W)=\bigoplus S_\lambda W,\qquad \bigwedge{} ^2(S^n W)=\bigoplus S_{\delta} W,$$
where the first sum runs over representations corresponding to $\lambda$ of weight $2n$ with two rows of even length and the second sum
runs over representations corresponding to $\delta$ of weight $2n$ with two rows of odd length.
\kob
\dow
Consider the multiplicity inside $S^2(S^n W)$ of a component with rows $\lambda_1, \lambda_2$, where $\lambda_1+\lambda_2=2n$ .
From Lemma \ref{Redlemsym} we know that this multiplicity is equal to the multiplicity of the component with one row of
length $\lambda_1-\lambda_2$ inside $\bigwedge^2(S^{n-\lambda_2} W)$ for $\lambda_2$ odd, and inside $S^2(S^{n-\lambda_2} W)$
for $\lambda_2$ even.
Hence it is one for $\lambda_2$ even and zero for $\lambda_2$ odd. A similar argument leads to the second equality.
\kdow
We proceed to the proof of Theorem \ref{th:plethysm3}. Due to the dualities \ref{dualities} we may consider only $S^\mu(S^k W)$ for $S^\mu$
a third symmetric or skew-symmetric power. Let us introduce some notation for symmetric polynomials.
\dfi[$h_k(x^a)$, $\psi_\alpha(h_k)$]
Consider $d$ variables $x_1,\dots,x_d$.
For $a\in\N$ let $h_k(x^a)$ be the complete symmetric polynomial of degree $k$ in the variables $x_1^a,\dots,x_d^a$.
We also define for a multi-index $\alpha$ of length $j$:
$$\psi_{\alpha}(h_k):=\prod_{i=1}^j h_k(x^{\alpha_{i}}).$$
\kdfi
The character of the representation $S^\mu(S^k W)$ equals $\sum \pm\frac{z_\alpha}{d!}\psi_\alpha(h_k)$, where the sum is taken over all
partitions $\alpha$ of $d$ and $z_\alpha$ is the number of permutations of combinatorial type $\alpha$ in the group $S_d$.
Our aim is to decompose $\psi_\alpha(h_k)$ into a sum of Schur polynomials. To do this we multiply $\psi_\alpha(h_k)$ by the discriminant
$\prod_{i<j} (x_i-x_j)$. Assume that there are $d-1$ variables. The coefficient of $s_\lambda$ inside $\psi_\alpha(h_k)$ equals the coefficient
of the monomial $x_1^{\lambda_1+d-2}\cdots x_{d-1}^{\lambda_{d-1}}$ in $\psi_\alpha(h_k)\prod_{i<j} (x_i-x_j)$ \cite[Appendix]{FultonRT}, \cite{macdonald1998symmetric}.

Let  $d=3$. By Lemma \ref{Redlemsym} we can assume that $\lambda$
has two rows $\lambda_1, \lambda_2$ with $\lambda_1+\lambda_2=3k$.
There are 3 partitions of the number $3$ to consider.
\begin{enumerate}
\item $3=1+1+1$. Here we need to compute the contribution of $h_3(x)^3$. This follows from Pieri's rule.
This contribution is equal to
$\lambda_2+1$ for $\lambda_2\leq k$, and $\lambda_1-\lambda_2+1$ for $\lambda_2\geq k$.
\item $3=2+1$. Note that the coefficient of a monomial $x_1^{3k-a}x_2^a$ in $h_k(x^2)h_k(x)$ for $a\leq 2k$ equals the number of
even integers less or equal to $a$ and greater or equal to $\max(0,a-k)$. We can easily deduce the coefficient of $x_1^{\lambda_1+1}x_2^{\lambda_2}$
in $(x_1-x_2)h_k(x^2)h_k(x)$. When $\lambda_2\leq k$, we get
$0$ for $\lambda_2$ odd, and $1$ for $\lambda_2$ even.
When $\lambda_2\geq k$, we get $0$ for $k$ odd, while for $k$ even we get $1$ if $\lambda_2$ is even, and $-1$ if $\lambda_2$ is odd.
\item $3=3$. The coefficient of $x_1^{3k-a}x_2^a$ in $h_k(x^3)$ is equal to $1$ if $a$ is divisible by $3$ and $0$ otherwise. Thus the
coefficient of $x_1^{\lambda_1+1}x_2^{\lambda_2}$ in $(x_1-x_2)S^k(x^3)$ is equal to $1$ if $\lambda_2=0 \;(\mathrm{mod}\;3)$,
$-1$ if $\lambda_2=1 \;(\mathrm{mod}\;3)$, $0$ if $\lambda_2=2 \;(\mathrm{mod}\;3)$.
\end{enumerate}
Finally, recall that the contribution form 1) is taken with coefficient $\frac{1}{6}$, from 2) with $\frac{1}{2}$ and from 3) with $\frac{1}{3}$.
This finishes the proof.
\kwadrat

\bibliographystyle{amsalpha}
\bibliography{Xbib}

\def\cprime{$'$} \def\cprime{$'$}
\providecommand{\bysame}{\leavevmode\hbox to3em{\hrulefill}\thinspace}
\providecommand{\MR}{\relax\ifhmode\unskip\space\fi MR }
\providecommand{\MRhref}[2]{%
  \href{http://www.ams.org/mathscinet-getitem?mr=#1}{#2}
}
\providecommand{\href}[2]{#2}
\begin{thebibliography}{MOZ12}

\bibitem[BCV13]{bruns2013relations}
Winfried Bruns, Aldo Conca, and Matteo Varbaro, \emph{Relations between the
  minors of a generic matrix}, Advances in Mathematics \textbf{244} (2013),
  171--206.

\bibitem[CGR]{plet2}
Y.M. Chen, A.M. Garsia, and J.~Remmel, \emph{Algorithms for plethysm, in
  combinatorics and algebra (boulder, colo., 1983)}, Contemp. Math.
  \textbf{34}, 109--153.

\bibitem[CT92]{Carre}
Christophe Carr{\'e} and Jean-Yves Thibon, \emph{Plethysm and vertex
  operators}, Adv. in Appl. Math. \textbf{13} (1992), no.~4, 390--403.

\bibitem[Dix77]{dixmier1977enveloping}
Jacques Dixmier, \emph{Enveloping algebras}, vol.~14, Newnes, 1977.

\bibitem[DK11]{DK2}
Jan Draisma and Jochen Kuttler, \emph{Bounded-rank tensors are defined in
  bounded degree}, arXiv:1103.5336v2 (2011).

\bibitem[Dra13]{draisma2013noetherianity}
Jan Draisma, \emph{Noetherianity up to symmetry}, arXiv preprint
  arXiv:1310.1705 (2013).

\bibitem[FH91]{FultonRT}
William Fulton and Joe Harris, \emph{Representation theory}, Graduate Texts in
  Mathematics, vol. 129, Springer-Verlag, New York, 1991, A first course,
  Readings in Mathematics.

\bibitem[FZ00]{fomin2000recognizing}
Sergey Fomin and Andrei Zelevinsky, \emph{Recognizing schubert cells}, Journal
  of Algebraic Combinatorics \textbf{12} (2000), no.~1, 37--57.

\bibitem[GS87]{gel1987combinatorial}
Izrail~Moiseevich Gel'fand and Vera~V Serganova, \emph{Combinatorial geometries
  and torus strata on homogeneous compact manifolds}, Russian Mathematical
  Surveys \textbf{42} (1987), no.~2, 133--168.

\bibitem[Lan12]{landsberg2012tensors}
Joseph~M Landsberg, \emph{Tensors:: Geometry and applications}, vol. 128, AMS
  Bookstore, 2012.

\bibitem[Lit44]{littlewood1944invariant}
Dudley~E Littlewood, \emph{On invariant theory under restricted groups},
  Philosophical Transactions of the Royal Society of London. Series A,
  Mathematical and Physical Sciences \textbf{239} (1944), no.~809, 387--417.

\bibitem[LM04]{MR2097214}
J.~M. Landsberg and L.~Manivel, \emph{On the ideals of secant varieties of
  {S}egre varieties}, Found. Comput. Math. \textbf{4} (2004), no.~4, 397--422.

\bibitem[LO11]{OttavianiLandsberg}
Joseph~M Landsberg and Giorgio Ottaviani, \emph{Equations for secant varieties
  of veronese and other varieties}, Annali di Matematica Pura ed Applicata
  (2011), 1--38.

\bibitem[LW07]{MR2363430}
J.~M. Landsberg and Jerzy Weyman, \emph{On tangential varieties of rational
  homogeneous varieties}, J. Lond. Math. Soc. (2) \textbf{76} (2007), no.~2,
  513--530.

\bibitem[LW09]{MR2533306}
\bysame, \emph{On secant varieties of compact {H}ermitian symmetric spaces}, J.
  Pure Appl. Algebra \textbf{213} (2009), no.~11, 2075--2086.

\bibitem[Mac98]{macdonald1998symmetric}
Ian~G Macdonald, \emph{Symmetric functions and hall polynomials}, Oxford
  University Press on Demand, 1998.

\bibitem[Man98]{MR1651092}
Laurent Manivel, \emph{Gaussian maps and plethysm}, Algebraic geometry
  ({C}atania, 1993/{B}arcelona, 1994), Lecture Notes in Pure and Appl. Math.,
  vol. 200, Dekker, New York, 1998, pp.~91--117.

\bibitem[Man09]{manivel2009spinor}
\bysame, \emph{On spinor varieties and their secants}, SIGMA \textbf{5} (2009),
  078.

\bibitem[MM12]{JaLaurent}
Laurent Manivel and Mateusz Micha{\l}ek, \emph{Effective constructions in
  plethysms and weintraub’s conjecture}, Algebras and Representation Theory
  (2012), 1--11.

\bibitem[MOZ12]{MOZ}
Mateusz Michalek, Luke Oeding, and Piotr Zwiernik, \emph{Secant cumulants and
  toric geometry}, arXiv preprint arXiv:1212.1515 (2012).

\bibitem[OR11]{RaicuOeding}
Luke Oeding and Claudiu Raicu, \emph{Tangential varieties of segre varieties},
  arXiv preprint arXiv:1111.6202 (2011).

\bibitem[Rai12]{Raicu}
Claudiu Raicu, \emph{Secant varieties of segre--veronese varieties}, Algebra \&
  Number Theory \textbf{6} (2012).

\bibitem[SZ11]{ZwiernikSturmfels}
Bernd Sturmfels and Piotr Zwiernik, \emph{Binary cumulant varieties}, Annals of
  Combinatorics (2011), 1--22.

\bibitem[Thr]{plet1}
R.~M. Thrall, \emph{On symmetrized kronecker powers and the structure of the
  free lie ring}, American Journal of Mathematics \textbf{64}, 371--388.

\bibitem[Zak05]{Zak}
FL~Zak, \emph{Tangents and secants of algebraic varieties}, vol. 127, AMS
  Bookstore, 2005.

\end{thebibliography}

\end{document}